\documentclass[12pt,dvips,reqno]{amsart}
\usepackage{euler, amsfonts, amssymb, latexsym, epsfig,epic,stmaryrd}

\newtheorem{Theorem}{Theorem} 
\newtheorem{Proposition}{Proposition} 
\newtheorem{Lemma}{Lemma}

\newtheorem{Corollary}{Corollary}
\newtheorem*{Corollary*}{Corollary}
 
\newtheorem*{Theorem*}{Theorem}
\theoremstyle{remark}
\newtheorem{Example}{Example}

\newcommand\iso{{\mathrel{\ \cong\ }}}

\newcommand\onto{\mathop{\twoheadrightarrow}}

\newcommand\into{\operatorname*{\hookrightarrow}}

\newcommand\nulset{\emptyset}

\theoremstyle{plain}

\renewenvironment{quotation}
{\list{}{
    \setlength\itemindent{0em}%
    \setlength\leftmargin{1.5em}
    \setlength\rightmargin{1.5em}
  }%
\item[]}
{\endlist}

\newcommand\dfn{\bf} 

\newcommand\Gm{{\mathbb G}_m}

\newcommand\kk{{\mathbb F}}

\newcommand\comment[1]{{\bf *** #1 ***}}
\newcommand\junk[1]{}

\renewcommand\AA{{\mathbb A}}

\begin{document}
\pagestyle{plain}

\title{Schubert patches degenerate to subword complexes}

\author{Allen Knutson}
\thanks{Supported by an NSF grant.}
\email{allenk@math.ucsd.edu}
\dedicatory{Dedicated to Bert Kostant on the occasion of his 80th birthday}
\date{Hallowe'en 2007}

\maketitle

\begin{abstract}
  We study the intersections of general Schubert varieties $X_w$
  with permuted big cells, and give an inductive degeneration of each
  such ``Schubert patch'' to a Stanley-Reisner scheme.
  Similar results had been known for Schubert patches in various
  types of Grassmannians. We maintain reducedness using the 
  results of [Knutson 2007] on automatically reduced degenerations,
  or through more standard cohomology-vanishing arguments.
  
  The underlying simplicial complex of the Stanley-Reisner scheme is a
  {\em subword complex}, as introduced for slightly different purposes
  in [Knutson-Miller 2004], and is homeomorphic to a ball.  
  This gives a new proof of the
  Andersen-Jantzen-Soergel/Billey and Graham/Willems formulae for
  restrictions of equivariant Schubert classes to fixed points.
\end{abstract}

{\small \tableofcontents}


Fix a pinning $(G,T,W,N_\pm,B_\pm=TN_\pm)$ 
of a complex reductive Lie group $G$.
Then there are correspondences between the Weyl group $W:=N(T)/T$,
the $B_-$-orbits on the flag manifold $G/B_+$,
and the set of $T$-fixed points $(G/B_+)^T$,
namely $w \mapsto B_- w B_+/B_+, wB_+/B_+$.
Let $X_w := \overline{B_- w B_+}/B_+ \subseteq G/B_+$ be
the {\dfn Schubert variety} associated to $w\in W$.

\subsection{Cohomology: Schubert classes in equivariant $K$-theory}

The $T$-invariant cycles $\{X_w\}$ define bases in many (co)homology theories
of $G/B$, in particular $T$-equivariant ones. We will focus our attention 
on equivariant $K$-theory $K_T(G/B)$, an algebra over 
the (Laurent polynomial) ring $K_T(pt)$ of virtual characters of the torus $T$.
However, the results we mention have (and imply) analogues 
for equivariant cohomology.

The structure sheaf of $X_w$ defines an element $S_w \in K_T(G/B)$.
These classes $S_w$ were first calculated by Demazure (implicitly, in
his formula \cite{Demazure} for the characters of Demazure modules),
though his calculation had a gap (later filled) which we will come to
in a moment.

In the paper \cite{KK} by Bertram Kostant and Shrawan Kumar, they 
suggest that one describe $S_w$ using the restriction map
(and $K_T(pt)$-algebra homomorphism)
$$ K_T(G/B) \into K_T\left((G/B)^T\right) = \bigoplus_{w\in W} K_T(pt) $$
The key point is that this map is {\em injective}, so no information
is lost by localizing to fixed points. To match notation with \cite{KK}
we think of the weight lattice $T^*$ additively, 
and denote the class of the $1$-dimensional representation with
weight $\lambda \in T^*$ by $e^\lambda \in K_T(pt)$.

Let $S_w|_v \in K_T(pt)$ denote the restriction of the class $S_w$ to the 
$T$-fixed point $v$.
Then $S_w|_v$ can be calculated inductively in $v$:

\begin{Theorem}
\label{thm:KK}
  Let $v,w$ be elements of $W$.
  If $S_w|_v \neq 0$, then $v\geq w$ in the Bruhat order.
  Assume this hereafter.

  If $v=1$, then $w=1$ and $S_w|_v = 1$.
  Otherwise, there exists a simple root $\alpha$ such that $v r_\alpha < v$
  in the Bruhat order (equivalently, $v\cdot\alpha$ is a negative root). 
  Then $S_w|_v$ can be computed 
  from $S_w|_{v r_\alpha}$ and $S_{w r_\alpha}|_{v r_\alpha}$:
  \begin{enumerate}
  \item If $w r_\alpha > w$, then $S_w|_v = S_w|_{v r_\alpha}$.
  \item If $w r_\alpha < w$ but $w \not\leq v r_\alpha$, then
    $ S_w|_v = (1 - e^{v \cdot \alpha}) S_{w r_\alpha}|_{v r_\alpha}. $
  \item If $w r_\alpha < w \leq v r_\alpha$, then
    $ S_w|_v = S_w|_{v r_\alpha} 
    + (1 - e^{v \cdot \alpha}) S_{w r_\alpha}|_{v r_\alpha}
    - (1-e^{v\cdot \alpha}) S_w|_{v r_\alpha}. $
  \end{enumerate}
  (In fact (3) includes (2).)
\end{Theorem}

This is essentially in \cite{KK} stated using Demazure operators,
and also follows trivially from the Graham/Willems formulae for $\{S_w|_v\}$
\cite{Graham,Willems}.
One can give a reasonably straightforward direct proof 
(in the finite-dimensional case) if one computes with the
{\em Bott-Samelson-Demazure-Hansen} resolutions of Schubert varieties.
However, such a proof depends on the fact that the map from the Bott-Samelson 
manifold to the Schubert variety takes the fundamental $K$-class to the 
fundamental $K$-class, which can be proven by showing that the
higher direct images vanish of the pushforward of the
structure sheaf (i.e. that Schubert varieties have rational singularities).
This was exactly the gap in Demazure's argument; an
account of this (and how the gap was eventually filled using
characteristic $p$ methods) can be found in \cite[chapter 14]{Jantzen}.

In theorem \ref{thm:degenlemmas} we will give a geometric interpretation
(and new proof)
of theorem \ref{thm:KK}, based on the degenerations of \cite{reduced}
applied to ``Schubert patches''.
In theorem \ref{thm:subword} we will give
a new combinatorial interpretation of theorem \ref{thm:KK}, based on 
vertex decompositions of the ``subword complexes'' of \cite{KM,subword}.
In our principal new theorem, theorem \ref{thm:main},
we will relate the geometry and the combinatorics directly by giving
a degeneration of Schubert patches to Stanley-Reisner schemes of
subword complexes.

\subsection{Geometry: Schubert patches}

One thing we do in this paper is to replace the use of resolutions 
with degenerations, 
and replace the cohomology-vanishing arguments with lemmas from 
\cite{reduced} about {\em automatically reduced} degenerations.
These lemmas are in turn based on the theory of {\em branchvarieties} 
\cite{AK}, though we will not need to inspect those underpinnings
in the present paper.

We define $X_w|_v := X_w \cap (v N_- B_+/B_+)$ as the intersection
of $X_w$ with the {\dfn permuted big cell} $v N_-B_+/B_+$, and call it a 
{\dfn Schubert patch} on $X_w$, as the $\{X_w|_v, v\geq w\}$ form an
affine open cover of the Schubert variety $X_w$. 
Each Schubert patch carries an action of $T$,
and the $T$-equivariant transverse pullback diagram
$$
\begin{array}{ccccc}
  &&X_w|_v&\into& X_w \\
  &&\downarrow&&\downarrow \\
  vB_+/B_+ &\into&v N_-B_+/B_+&\into& G/B_+ 
\end{array}
$$
shows that we can compute the restriction $S_w|_v$ 
as the class $[X_w|_v] \in K_T(vN_-B_+/B_+) \iso K_T(pt)$. 
This is convenient for a number of purposes, one being that
the permuted big cell $vN_-B_+/B_+$ is a vector space whose $T$-weights
all live in the interior of a half-space of $T^*$. As such, $X_w|_v$ has a 
multigraded Hilbert series, and the $K_T$-class $S_w|_v$ is essentially 
this series times $v\cdot$(the Weyl denominator); see \cite[chapter 8.2]{MS}. 
This is the viewpoint of \cite{GhoRag,KodRag,KreLak,RU1,RU2}.

The degenerations we use are of a very specific type, which we
christened {\dfn geometric vertex decompositions} in \cite{gvd}.
As explained in \cite[section 4.1]{reduced},
the permuted big cell $v N_- B_+/B_+ \iso N_-$ factors $T$-equivariantly 
as a product of a line $L$ with weight $-v\cdot \alpha$ and 
a complementary hyperplane $H$. Let the multiplicative group $\Gm$
act on the permuted big cell $H\times L$ by scaling $L$, and define
$$ X' := \lim_{t\to 0}\ t\cdot X_w|_v $$
as the scheme-theoretic limit.
Automatically, $X'$ has the same multigraded Hilbert series and
(equivalently) defines the same $K_T$-class as $X_w|_v$.

We can now state our degeneration-based analogue of theorem \ref{thm:KK}.
In it, we use the notation $\AA^1_\lambda$ to denote the $1$-dimensional
$T$-representation with weight $\lambda$, and all isomorphisms
stated are $T$-equivariant.

\begin{Theorem}
\label{thm:degenlemmas}
  Let $v,w$ be elements of $W$.
  If $X_w|_v \neq \nulset$, then $v\geq w$ in the Bruhat order. 
  Assume this hereafter.

  If $v=1$, then $w=1$ and $X_w|_v = N_- B_+/B_+$. 
  Otherwise, there exists a simple root $\alpha$ such that $v r_\alpha < v$
  in the Bruhat order. 
  Let $X'$ be the degeneration of $X_w|_v$ described above.
  \begin{enumerate}
  \item If $w r_\alpha > w$, then $X' = X_w|_v$ (the limiting process
    is trivial), and 
    $$ X_w|_v \iso \Pi \times \AA^1_{-v\cdot\alpha}, \quad
    X_w|_{v r_\alpha} \iso \Pi \times \AA^1_{v\cdot\alpha} $$
    for the same $\Pi$.
  \item If $w r_\alpha < w$ but $w \not\leq v r_\alpha$, then
    $X' = X_w|_v$ (again, the limiting process is trivial), and 
    $$ X_{w r_\alpha}|_{v r_\alpha} \iso X_w|_v \times \AA^1_{v\cdot\alpha}.$$
  \item If $w r_\alpha < w \leq v r_\alpha$, then $X'$ is reduced,
    and has two components:
    $$ X' = \left(\Pi \times \{0\}\right) 
    \cup_{\Lambda \times \{0\}} 
    \left(\Lambda \times \AA^1_{-v\cdot\alpha} \right) $$
    where $\Pi \times \AA^1_{v\cdot \alpha} \iso X_{w r_\alpha}|_{v r_\alpha}$
    and $\Lambda \times \AA^1_{v\cdot \alpha} \iso X_w|_{v r_\alpha}$.
  \end{enumerate}
\end{Theorem}

\begin{proof}
  Parts (1) and (2) are proposition 6 of \cite{reduced} and are surely
  well-known to the experts. Part (3), which is much deeper, 
  is proposition 7 of \cite{reduced}. (This is why we separated cases
  (2) and (3) in theorem \ref{thm:KK}.)
\end{proof}

By manipulating Hilbert series, it is easy to recover each part of
theorem \ref{thm:KK} from the corresponding part of theorem 
\ref{thm:degenlemmas}. 
And indeed,
while the set-theoretic description of $X'$ is reasonably straightforward
\cite[theorem 2.2b]{gvd},
its reducedness is essentially equivalent to the Kostant-Kumar recursion.
But where part (3) of theorem \ref{thm:KK} was proven using the
vanishing of higher cohomology of the Bott-Samelson-Demazure-Hansen resolution,
part (3) of theorem \ref{thm:degenlemmas} is based on 
\cite[proposition 7]{reduced}, which shows that the (by definition
reduced) ``limit branchvariety'' coincides with the limit subscheme.

Having a degeneration implies more than merely an equality of $K_T$-classes: 
in \cite{reduced} we use these results to give a new proof that
Schubert varieties are normal and Cohen-Macaulay.

In the remainder of the paper it will be convenient to work not directly
with Schubert patches, but the more economical 
{\dfn Kazhdan-Lusztig varieties} $X_{w\circ}^v := X_w \cap X^v_\circ$,
where $X^v_\circ$ denotes the {\dfn opposite Schubert cell} $N_+ v B_+/B_+$.
These are obtained from Schubert patches by factoring out a largely
irrelevant vector space:

\begin{Lemma}\cite[Lemma A.4]{KazhdanLusztig}\label{lem:KLvspatches}
  Let $w,v \in W$. Then there is a $T$-equivariant factorization
  $$ X_w|_v \iso X_{w\circ}^v \times X_v^\circ $$
  where $X_v^\circ := N_- v B_+/B_+$ is just a $T$-vector space, with weights
  $\{ \beta<0 : v\cdot \beta < 0\}$. In particular the dimension
  of the Kazhdan-Lusztig variety is
  $ \dim X_{w\circ}^v = \ell(v) - \ell(w).$
\end{Lemma}

Each result about Schubert patches has an equivalent, though often simpler, 
version for Kazhdan-Lusztig varieties.
The only possibly subtle one is the restriction of
the geometric vertex decomposition, which is a family of subschemes
of $X_1|_v \cup X_1|_{v r_\alpha}$; 
for this one intersects each fiber of the family with $X^v$,
obtaining a family of subschemes of $X^v_\circ \cup X^{v r_\alpha}_\circ$.

\newtheorem*{thmprime}{Theorem \ref{thm:degenlemmas}'}
\begin{thmprime}
  Let $v,w$ be elements of $W$.
  If $X_{w\circ}^v \neq \nulset$, then $v\geq w$ in the Bruhat order. 
  Assume this hereafter.

  If $v=1$, then $w=1$ and $X_{w\circ}^v = B_+/B_+$. 
  Otherwise, there exists a simple root $\alpha$ such that $v r_\alpha < v$
  in the Bruhat order. 
  \begin{enumerate}
  \item If $w r_\alpha > w$, then 
    $ X_{w\circ}^v 
    \iso X_{w\circ}^{v r_\alpha} \times \AA^1_{-v\cdot\alpha}. $
  \item If $w r_\alpha < w$ but $w \not\leq v r_\alpha$, then
    $ X_{w\circ}^v \iso X_{w r_\alpha \circ}^{v r_\alpha}.$
  \item If $w r_\alpha < w \leq v r_\alpha$, then $X'$ is reduced,
    and has two components:
    $$ X' = \left(\Pi \times \{0\}\right) 
    \cup_{\Lambda \times \{0\}} 
    \left(\Lambda \times \AA^1_{-v\cdot\alpha} \right) $$
    where $\Pi \iso X_{w r_\alpha \circ}^{v r_\alpha}$
    and $\Lambda \iso X_{w\circ}^{v r_\alpha}$.
  \end{enumerate}
\end{thmprime}

\subsection{Combinatorics: subword complexes}

When attempting to unwind theorem \ref{thm:KK} to a direct formula for
$S_w|_v$, as in Kumar's appendix to \cite{Billey},
one is led naturally to the definition of a {\em subword complex}
\cite{subword} (though our motivation at the time was slightly different).

Let $Q = (\alpha_1,\alpha_2,\ldots,\alpha_k)$ be a sequence of simple roots
such that $v = \prod_{i=1}^k r_{\alpha_i}$, and $k$ is minimized.
Then $Q$ is called a {\dfn reduced word} for $v$,
and $k$ its {\dfn length}, usually denoted $\ell(v)$.
(Warning: because we use the Kostant-Kumar recurrence based on 
$v r_\alpha$ and not one based on $r_\alpha v$, 
the first root used in applying the recurrence is $\alpha_k$, not $\alpha_1$.)

In \cite{KM,subword} we defined the {\dfn subword complex} $\Delta(Q,w)$ 
associated to a reduced word $Q$ and a Weyl group element $w$ 
as the simplicial complex whose vertex set is $Q$ (or really, $1\ldots k$)
with $F \subseteq Q$ a {\dfn facet} (maximal face) iff the complement
$Q\setminus F$ is a reduced word for $w$.

Even when $Q\setminus F$ is not a reduced word, we can define its
{\dfn Demazure product} by multiplying the reflections in order, 
omitting along the way any one that brings us lower in the Bruhat order.
(Equivalently, one may take the product of any maximal reduced subword.)

\begin{Proposition}\cite{subword}\label{prop:subwordpaper}
  The subword complex $\Delta(Q,w)$ is homeomorphic to a ball or sphere;
  in particular every {\dfn ridge} (codimension $1$ face)
  is contained in one or two facets. 

  For any face $F\in \Delta(Q;w)$, the Demazure product of $Q\setminus F$
  is $\geq w$ in the Bruhat order, with equality iff $F$ is an
  {\dfn interior} face (i.e. if $F$ is contained in no ridge contained in only
  one facet). 
\end{Proposition}

To any simplicial complex $\Delta$ on vertex set $Q$, and a field $\kk$,
one may associate the (affine) {\dfn Stanley-Reisner scheme} 
$SR(\Delta) \subseteq \kk^Q$, the 
union of the corresponding coordinate planes:
$$ SR(\Delta) := \bigcup_{S \in \Delta} \kk^S. $$
These schemes are invariant under the action of the torus $(\Gm)^Q$ that
dilates the coordinates independently. (Indeed, they are characterized
by this invariance plus their reducedness; note too that $\Delta$ can
be reconstructed from $SR(\Delta)$.)
As such $SR(\Delta)$ has an associated 
multigraded Hilbert series in the variables $(q_1,\ldots,q_{|Q|})$:
$$ h_{SR(\Delta)} 
= \sum (\hbox{those monomials whose variables form a face in $\Delta$}) 
= \sum_{F\in \Delta} \prod_{q_j\in F} \frac{q_j}{1-q_j}. $$
Equivalently, one can compute the class 
$k_\Delta := [SR(\Delta)] \in K_{(\Gm)^Q}(\kk^Q)$, by
$$ k_\Delta 
= \left(\prod_{q_j} (1-q_j)\right) h_{SR(\Delta)}
= \sum_{F\in \Delta} \prod_{q_j\in F} q_j \prod_{q_j\notin F} (1-q_j). $$
One can give an alternate formula for the Hilbert series $h_{SR(\Delta)}$ 
in which the summands are of the form $\prod_{q_j\in F} 1/(1-q_j)$
rather than $\prod_{q_j\in F} q_j/(1-q_j)$, corresponding to writing
$\Delta$ as a union of closed faces rather than open faces.  
The resulting inclusion-exclusion of the faces is particularly simple
in the case of $\Delta$ a ball or sphere, and becomes an alternating
sum over the interior faces:

\junk{
  A {\dfn cone point} of a simplicial complex $\Delta$ is a vertex that
  occurs in every facet, and as such doesn't change the combinatorics in
  any interesting way.
  For technical reasons, in this paper we will use a tiny variation
  of the definition of subword complex that introduces many cone points.
  Let $Q$ be a reduced word for $w_0$. Then let $\Delta(Q,i,w)$
  be the complex with vertex set $Q$, and $F\subseteq Q$ a facet
  if the complement {\em uses only the first $i$ letters} and is
  a reduced word for $w$. Put another way, $\Delta(Q,i,w)$ is a multicone
  on the usual subword complex $\Delta($first $i$ letters of $Q,w)$,
  with the last $\ell(w_0)-i$ vertices attached as cone points. One effect
  of this is that for any facet $F$, $|F| = \dim X_w$, independent of $Q$.
}

\begin{Corollary}\cite[Lemma 4.2]{subword}\label{cor:interior}
  $$ h_{SR(\Delta(Q,w))}
  = \sum_{F \subseteq Q \atop \prod (Q\setminus F) = w} 
  (-1)^{|Q\setminus F|} \prod_{q_j\in F} \frac{1}{1-q_j} $$
  where $\prod (Q\setminus F)$ means the Demazure product of the subword.
\end{Corollary}

We comment that if one sums over {\em all} faces of $\Delta(Q,w)$, not
just the interior ones, one obtains the $K_T$-class of the ideal sheaf
of the ``boundary'' of the Schubert variety. These ideal sheaves give
another $K_T(pt)$-basis of $K_T(G/B)$ which (upon twisting by the
long element in $W$)
are the dual basis under the Poincar\'e pairing on $K_T(G/B)$,
and this formula for them appears in \cite[section 3]{Graham}.

To give a simplicial-complex analogue of theorem \ref{thm:KK}, 
we will need an analogue of the decomposition that appears in its part (3).
The {\dfn deletion} of a vertex $p$ from a simplicial complex $\Delta$ 
is the subcomplex $del_p \Delta := \{F\in \Delta : F\not\ni p\}$,
and the {\dfn star} of the vertex $p$ is the subcomplex 
$star_p \Delta := \{F \in \Delta : F\cup \{p\} \in \Delta\}$.  
So $p$ is a cone point of its star, and deleting it we get the {\dfn link} 
$link_p \Delta$ of $p$. The decomposition 
$$ \Delta = del_p \Delta \cup_{link_p \Delta} star_p \Delta $$
is a {\dfn vertex decomposition}, as used in \cite{BilleraProvan}.
At the risk of getting ahead of ourselves, we mention that
this should be seen as analogous
to part (3) of theorem \ref{thm:degenlemmas}
(hence the term ``geometric vertex decomposition''), with 
$del_p \Delta$ and $link_p \Delta$ playing the roles of $\Pi$ and $\Lambda$.

\begin{Theorem}\cite{subword}\label{thm:subword}
  Let $v,w$ be elements of $W$. Let $Q$ be a reduced word for $v$. 
  If $\Delta(Q,w) \neq \nulset$, then $v\geq w$ in the Bruhat order. 
  Assume this hereafter.

  If $v=1$, then $w=1$, $Q = ()$ and $\Delta(Q,w) = \{\nulset\}$.
  Otherwise, let $\alpha$ be the last simple root listed in $Q$,
  and let $Q'$ be $Q$ with this root dropped. (Hence $v r_\alpha < v$.)
  \begin{enumerate}
  \item If $w r_\alpha > w$, then $\ell(v)$ is a cone point of 
    $\Delta(Q,w)$, and 
    $$ link_{\ell(v)} \Delta(Q,w) = del_{\ell(v)} \Delta(Q,w) = \Delta(Q',w).$$
  \item If $w r_\alpha < w$ but $w \not\leq v r_\alpha$, 
    then no face of $\Delta(Q,w)$ uses $\ell(v)$, and
    $$ \Delta(Q,w) = del_{\ell(v)} \Delta(Q,w) = \Delta(Q',w r_\alpha). $$
  \item If $w r_\alpha < w \leq v r_\alpha$, then the vertex decomposition
    at the vertex $\ell(v)$ is into
    $$     link_{\ell(v)} \Delta(Q,w) = \Delta(Q',w), \qquad
    del_{\ell(v)} \Delta(Q,w) =   \Delta(Q',w r_\alpha).
    $$
  \end{enumerate}
\end{Theorem}

\begin{proof}
  The proofs are all largely tautologies based on the Bruhat order,
  taking care not to be confused by the {\em complementation} involved
  in the definition of subword complex.
  \begin{enumerate}
  \item If $w r_\alpha > w$, then the last letter in a reduced word for $w$
    cannot be $r_\alpha$. Complementing, $\ell(v)$ must lie in every facet
    of $\Delta(Q,w)$. Then any subword of $Q$ using the
    first $\ell(v)$ letters but avoiding $\ell(v)$ is equivalently a
    subword of $Q$ using the first $\ell(v)-1$ letters, i.e. of $Q'$.
  \item If $w \not\leq v r_\alpha$, then no subword of the first $\ell(v)-1$
    letters in $Q$ has product $w$. Hence any subword of the first $\ell(v)$
    letters with product $w$ must use the $\ell(v)$th letter.
    Removing that letter, we get the product $w r_\alpha$.
  \item As the answer suggests, this is essentially a combination of
    the previous two arguments.
  \end{enumerate}
  (As in theorem \ref{thm:KK}, (2) is really a subproblem of (3).
  It is only in theorem \ref{thm:degenlemmas} that it is in any way
  natural to separate them.)  
\end{proof}

We are ready to connect theorems \ref{thm:KK} and \ref{thm:subword}:
the equivariant $K$-classes computed in theorem \ref{thm:KK} can
be computed from the Hilbert series associated to subword complexes.

\begin{Corollary}\label{cor:hilbmatch}
  Fix $Q,v,w$ as above, and let 
  $$ \beta_j := \left( \prod_{k=1}^{j-1} r_{\alpha_k} \right)\cdot \alpha_j,
  \qquad j=1,\ldots,\ell(v). $$
  Then the specialization
  $$ h_{Q,w} := h_{SR(\Delta(Q,w))}\hbox{ with } q_i \mapsto e^{\beta_i} $$
  of the Hilbert series of the Stanley-Reisner scheme of the subword
  complex is the $T$-equivariant Hilbert series of the Kazhdan-Lusztig
  variety $S_{w\circ}^v$.
  Equivalently, 
  $$ S_w|_v = \left(\prod_{\beta > 0 \atop v\cdot \beta>0} 
    (1-e^{-v\cdot \beta})\right) h_{Q,w} $$
  where the product is over those positive roots $\beta$ of $G$ that
  stay positive when twisted by $v$.
\end{Corollary}

\begin{proof}
  Assume $v\geq w$, for otherwise both sides are zero. 
  Then if $v=1$, both sides are $1$. Otherwise $\ell(v)\geq 1$
  and as before we let $Q'$ be $Q$ minus its last letter.

  Theorem \ref{thm:subword} implies corresponding results for the
  Hilbert series and $K$-polynomials:
  \begin{enumerate}
  \item If $w r_\alpha > w$, then
    $$ h_{\Delta(Q,w)} = h_{\Delta(Q',w)} / (1 - q_{\ell(v)}), \qquad
    k_{Q,w} = k_{Q',w}. $$
  \item If $w r_\alpha < w$ but $w \not\leq v r_\alpha$, then
    $$ h_{\Delta(Q,w)} =  h_{\Delta(Q',w r_\alpha)}, \qquad
    k_{Q,w} = (1 - q_{\ell(v)}) k_{Q',w r_\alpha}. $$
  \item If $w r_\alpha < w \leq v r_\alpha$, then
    $$ h_{\Delta(Q,w)} 
    = h_{\Delta(Q',w)} / (1 - q_{\ell(v)}) + h_{\Delta(Q',w r_\alpha)} 
    -  h_{\Delta(Q',w)} $$
    $$ k_{Q,w} = k_{Q',w} + (1 - q_{\ell(v)}) k_{Q',w r_\alpha} 
    - (1 - q_{\ell(v)}) k_{Q',w}. $$
  \end{enumerate}
  Under the specialization $q_i \mapsto e^{\beta_i}$ of $k_{Q,w}$,
  we recover the equations from theorem \ref{thm:KK}. 
  So under this specialization, $S_w|_v = k_{Q,w}$,
  and $k_{Q,w} = \left(\prod_{\beta > 0 \atop v\cdot \beta>0} 
    (1-e^{-v\cdot \beta})\right) h_{Q,w}$.
  This establishes the second claim.
  
  For the first, we use the transversality of the pullback diagram
  $$
  \begin{array}{ccc}
    X_{w\circ}^v&\into& X_w|_v \\
    \downarrow&&\downarrow \\
    N_+ v B_+/B_+&\into& v N_-B_+/B_+
  \end{array}
  $$
  and the relation of Hilbert series to $K_T$-classes \cite[chapter 8.2]{MS},
  we can compute 
  $$ S_w|_v 
  = h_{X_w|_v} \prod_{\beta < 0} (1 - e^{v\cdot \beta})
  = h_{X_{w\circ}^v} \prod_{\beta < 0, v\cdot \beta < 0} (1 - e^{v\cdot \beta})
  $$
  where the products account for the weights on the $T$-spaces
  $N_+ v B_+/B_+$, $v N_-B_+/B_+$ respectively.
\end{proof}

\begin{Corollary}\cite{Graham,Willems}
  Let $Q$ be a reduced word for $v$, and $\Delta(Q,w)$ the subword complex.
  Let $\Delta(Q,w)^\circ$ be the set of interior faces, i.e. those $F$ such 
  that the Demazure product of $Q\setminus F$ is exactly $w$ (not $>w$). Then
  $$ S_w|_v 
  = \sum_{F \in \Delta(Q,w)^\circ} (-1)^{|Q\setminus F|}
  \left( \prod_{i=1}^{\ell(v)} \widehat{(1 - e^{-\alpha_i})}^{[i\in F]}
    \ r_{\alpha_i}\right) \cdot 1 $$
  where $\widehat x$ is the multiply-by-$x$ operator, and
  $[i\in F]=0,1$ according to whether the condition fails or is satisfied.
\end{Corollary}

\begin{proof}
  This is the combination of corollary \ref{cor:interior} and the
  second half of corollary \ref{cor:hilbmatch}.
\end{proof}

This $K$-theory result implies in turn the corresponding equivariant
cohomology result from \cite[Appendix D]{AJS} and \cite{Billey}, which is a 
sum only over the facets of $\Delta(Q,w)$, rather than all interior faces.

\junk{
\begin{Corollary}
  Let $G/P$ be a {\dfn cominuscule} flag manifold, meaning that there
  exists a simple coroot $\vec v$ such that 

  \comment{should this go in this paper? minor literature search}
\end{Corollary}
}

The following technical lemma gives an inductive way to construct
subword complexes. Since it is purely combinatorial, we put it in
this section, but its main use will be in the geometry of the next section.

\begin{Lemma}\label{lem:inductivesubword}
  Fix $w\leq v \in W$, and let $Q$ be a reduced word for $v$.
  
  For each $i \leq \ell(v)$, let $v_i$ denote the Demazure product of
  the initial subword $(q_1,\ldots,q_i)$, and 
  $$ C_i 
  := \left\{\left(w' \leq v_i,\ S \subseteq \{i+1,\ldots,\ell(v) \}\right)\ :\ 
    \ell(w') + |S| = w,\ w' \prod S = w \right\}. $$
  By the length considerations, in each pair $(w',S) \in C_i$ the subword
  $S$ is automatically a reduced word for $w'^{-1} w$.
  Plainly $C_{\ell(v)} = \{ (w,\emptyset) \}$
  whereas $C_0 = \{ (1,S \subseteq Q) : S$ is a reduced word for $w\}$.

  There is a surjection $C_{i-1} \onto C_i$, taking
  $$ (w',S) \mapsto 
  \begin{cases}
    (w',S) & \hbox{if } i\notin S \\
    \left(w' r_\alpha, S\setminus \{i\}\right) & \hbox{if } i\in S
  \end{cases}
  $$

  Consequently, we can construct $C_{i-1}$ from $C_i$ as the
  disjoint union of the fibers of this surjection.
  Let $(w',S) \in C_i$, and let $\alpha$ be the $i$th root in the
  reduced word $Q$, so $v_i = v_{i-1} r_\alpha$.
  Then each fiber has one or two elements:
  \begin{enumerate}
  \item If $w' r_\alpha > w'$, then $w' \leq v_{i-1}$, and 
    $(w', S) \in C_{i-1}$ as well.
  \item If $w' r_\alpha < w'$ but $w' \not\leq v_i r_\alpha = v_{i-1}$, 
    then $(w' r_\alpha, \{i\} \cup S) \in C_{i-1}$.
  \item If $w' r_\alpha < w' \leq v_i r_\alpha = v_{i-1}$, 
    then $(w', S),(w' r_\alpha, \{i\} \cup S)$ are both in $C_{i-1}$.
  \end{enumerate}
\end{Lemma}

\begin{proof}
  All the claims made are essentially tautological.
\end{proof}

\subsection{From \#2 to \#3: Schubert patches degenerate to subword complexes}

We come shortly to the principal new theorem, theorem \ref{thm:main},
after an abbreviated history of related results.

There has been a great deal of work on degenerations of Schubert 
{\em varieties} (rather than patches) to unions of toric varieties,
starting with Hodge's degeneration of the Grassmannian (and its
Schubert varieties) in its Pl\"ucker embedding, limiting to what we
today would call the projective Stanley-Reisner scheme of the order
complex of the Bruhat order on the Grassmannian (see e.g. \cite{Hodge}).

For more general Schubert varieties in more general embeddings, 
it has been very fruitful to degenerate to unions not just of
projective spaces, but of more
complicated toric varieties, e.g. \cite{Chirivi,Caldero,KogM}.
(The degeneration in \cite{Chirivi} of a flag manifold {\em is} to a 
Stanley-Reisner scheme if the original flag manifold is a ``minuscule'' 
flag manifold in its fundamental embedding, but not otherwise.)
Unlike the geometric results in this paper, the constructions of these
degenerations for general $G$ have depended on deep algebraic results about
Lusztig's or Kashiwara's canonical bases.

Much less seems to be known if one insists on Stanley-Reisner schemes.
However one may change the game: rather than degenerating Schubert varieties,
one may degenerate matrix Schubert varieties \cite{KM,gvd} or
Schubert patches \cite{GhoRag,KodRag,KreLak,RU1,RU2}. One difference
when working with patches is that the choice of embedding becomes immaterial.

The Schubert patches considered to date in this context all live in various 
flavors of Grassmannians, i.e. minimal flag manifolds of classical groups.
We give now a uniform result producing Stanley-Reisner degenerations
of arbitrary Schubert patches in finite-dimensional $G/B$. Or rather,
we give a result about Kazhdan-Lusztig varieties, which one
may multiply by a vector space using lemma \ref{lem:KLvspatches}
if one prefers to work with Schubert patches.

\begin{Theorem}\label{thm:main}
  Fix $w\leq v \in W$, and let $Q$ be a reduced word for $v$.

  Then there is a sequence of flat $T$-equivariant degenerations, 
  starting from the Kazhdan-Lusztig variety $X_{w\circ}^v$, 
  and culminating in $SR(\Delta(Q,w))$.
\end{Theorem}

\begin{proof}
  The proof is of course inductive, and we need first to describe the
  structure of the intermediate cases, in a setting partway between
  Kazhdan-Lusztig varieties and Stanley-Reisner schemes of subword complexes.
  
  Let $C_i$ be the set of pairs $(w',S)$ defined in 
  lemma \ref{lem:inductivesubword}, and let
  $$ X_i := \bigcup_{(w',S) \in C_i} \left(
    X_{w' \circ}^{v_i} \times (\AA^1)^{ \{i+1,\ldots,\ell(v) \} \setminus S} 
  \right) \quad \subseteq X_\circ^{v_i} \times (\AA^1)^Q $$
  where $(\AA^1)^{ \{i+1,\ldots,\ell(v) \} \setminus S}$  
  denotes the evident coordinate subspace of $(\AA^1)^Q$. 
  Then $X_{\ell(v)} = X_w|_v \times \{0\}$ 
  and $X_0 = (B_+/B_+) \times SR(\Delta(Q,w))$.

  We can now give the correct inductive claim: 
  for each $i = 1,\ldots,\ell(v)$, 
  there is a flat $T$-equivariant degeneration of $X_i$ to $X_{i-1}$.
  Specifically, we will show that the degeneration described before 
  theorem \ref{thm:degenlemmas}, when applied to the first factor of $X_i$, 
  gives $X_{i-1}$.
  
  We first make a general comment about degenerating unions of closed
  subschemes (here the components of $C_i$).  Set-theoretically, the
  limit of a union is the union of the limits, but scheme-theoretically 
  there is usually only an inclusion. (Consider two points colliding
  in a line, whose scheme-theoretic limit is a fat point, containing
  the reduced union of the limit point with itself.)
  
  In the case at hand,
  we can follow an individual component $X_{w' \circ}^{v_i} \subseteq X_i$
  using theorem \ref{thm:degenlemmas}', and see that it produces
  exactly the components listed in lemma \ref{lem:inductivesubword}
  in the corresponding fiber of the $C_{i-1} \onto C_i$ surjection.
  By the above comment, we have shown that $X_i$ degenerates to a scheme whose
  reduction is $X_{i-1}$. 

  This gives an inequality on Hilbert series,
  with equality exactly if the degeneration is already reduced.
  Chaining these inequalities together, we get an inequality relating
  the Hilbert series of $X_0$ and $X_{\ell(v)}$. 
  But by corollary \ref{cor:hilbmatch}, we know these Hilbert series
  are equal. Hence each intermediate degeneration is indeed 
  of $X_i$ to $X_{i-1}$, scheme-theoretically.

  Chaining these degenerations $X_{\ell(v)} \leadsto \ldots \leadsto X_0$ 
  together, we have the sequence claimed in the theorem.
\end{proof}

Hartshorne's connectedness theorem for Hilbert schemes states that two
subschemes of the same projective space with the same $K$-class
can be connected by a series of deformations and degenerations.
So one might expect the Graham/Willems formula to imply the above theorem
directly. But this theorem is better in two ways: all the degenerations
preserve the {\em equivariant} $K$-class, and the theorem uses only 
degenerations (general to special), not deformations (special to general). 
So one may use semicontinuity arguments, e.g. the Stanley-Reisner
schemes of subword complexes being Cohen-Macaulay (since the
complexes are shellable) implies that Schubert patches are Cohen-Macaulay.
(In \cite{reduced} this argument was used one degeneration at a time.)

It seems likely (though we didn't pursue it) that one could use the
reduced word $Q$ to define good coordinates on the opposite cell $X_\circ^v$,
within which the above degeneration is by a Gr\"obner basis with
squarefree initial terms, and
that the reduced Gr\"obner basis could be inductively constructed
using the vertex decomposition of the subword complex.
A basis for this ideal (though not in these specific coordinates)
was already constructed in \cite[proposition 9.6.1]{LLM}, using
Frobenius splitting and canonical basis techniques.

\subsection{Acknowledgements}
Primarily I thank Rebecca Goldin, with whom we discovered the
Graham/Willems formula sometime between \cite{Graham} and \cite{Willems}.
This work would have been quite impossible without Ezra Miller, who
taught me so much about simplicial complexes. Thanks also to
Bill Graham for sending me his preprint \cite{Graham} and the
reference to \cite{AJS}.

\bibliographystyle{alpha}    

\end{document}